\documentclass[a4paper,reqno,10pt]{amsart}
\usepackage{newcent} 
\usepackage[T1]{fontenc}
\usepackage[centertags]{amsmath}
\usepackage{amsfonts}
\usepackage{amssymb}
\usepackage{amsthm}
\usepackage{newlfont}
\usepackage{graphicx}
\usepackage{color}
\usepackage{nogap-private}
\usepackage{mathtools}
\mathtoolsset{showonlyrefs}

\numberwithin{equation}{section}

\newtheorem{defin}{Definition}[section]
\newtheorem{theorem}[defin]{Theorem}
\newtheorem{lemma}[defin]{Lemma}

\theoremstyle{definition} {\newtheorem{remark}[defin]{Remark}}

\begin{document}

\title[Formulas for Relaxed Disarrangement Densities]{Explicit Formulas for Relaxed Disarrangement Densities Arising from Structured Deformations}
\author{Ana Cristina Barroso}
\address{Faculdade de Ci\^{e}ncias da Universidade de Lisboa, Departamento de Matem\'{a}tica and CMAF, Campo Grande, Edif\'{i}cio C6, Piso 1, 1749-016 Lisboa, Portugal}
\email[A.~C.~Barroso]{acbarroso@ciencias.ulisboa.pt}
\author{Jos\'{e} Matias}
\address{Departamento de Matem\'atica, Instituto Superior T\'ecnico, Av.\@ Rovisco Pais, 1, 1049-001 Lisboa, Portugal}
\email[J.~Matias]{jose.c.matias@tecnico.ulisboa.pt}
\author{Marco Morandotti}
\address{SISSA -- International School for Advanced Studies, Via Bonomea, 265, 34136 Trieste, Italy}
\email[M.~Morandotti \myenv]{marco.morandotti@sissa.it}
\author{David R.\@ Owen}
\address{Department of Mathematical Sciences, Carnegie Mellon University, 5000 Forbes Ave., Pittsburgh, PA 15213 USA}
\email[D.~R.~Owen]{do04@andrew.cmu.edu}
\date{\today} 

\maketitle

\begin{abstract}
In  this  paper  we  derive  explicit  formulas for disarrangement
densities of submacroscopic separations,
switches, and interpenetrations in the context of first-order structured
deformations.  Our derivation employs relaxation
within one mathematical setting for structured deformations of a specific,
purely interfacial density, and the formula we obtain agrees with one
obtained earlier in a different setting for structured deformations.
Coincidentally, our derivation provides an alternative method for
obtaining the earlier result, and we establish new explicit formulas for
other measures of disarrangements that are significant in applications.
\end{abstract}
\medskip

{\small
\keywords{\noindent {\bf Keywords:} {Structured deformations, relaxation, disarrangements, interfacial density, bulk density, isotropic vectors.}
}
\smallskip

\subjclass{\noindent {\bf {2010} 
Mathematics Subject Classification:}
{49J45 
(74A60, 
74G65, 
15A99). 
}}
}

\tableofcontents

\section{Introduction}

Structured deformations provide a multiscale geometry that
captures the contributions at the macrolevel of both smooth geometrical
changes and non-smooth geometrical changes (disarrangements) at
submacroscopic levels. \ For each (first-order) structured deformation $(g,G)$ 
of a continuous body, the tensor field $G$ is known to be a measure of
deformations without disarrangements, and $M:=\nabla g-G$ is known to be a
measure of deformations due to disarrangements. \ The tensor fields $G$ and 
$M$ together deliver not only standard notions of plastic deformation, but $M$
and its curl deliver the Burgers vector field associated with closed
curves in the body and the dislocation density field used in describing
geometrical changes in bodies with defects. Recently, Owen and Paroni
\cite{OwenParoni} evaluated explicitly some relaxed energy densities arising in
Choksi and Fonseca's energetics of structured deformations \cite{choksifonseca} and
thereby showed: (1) $\left( trM\right) ^{+}$, the positive part of $trM$, is
a volume density of disarrangements due to submacroscopic separations, 
(2) $\left( trM\right) ^{-}$, the negative part of $trM$, is a volume density of
disarrangements due to submacroscopic switches and interpenetrations, and
(3) $\left\vert trM\right\vert $, the absolute value of $trM$, is a volume
density of all three of these non-tangential disarrangements: separations,
switches, and interpenetrations. \ The main contribution of the present
research is to show that a different approach to the energetics of
structured deformations, that due to Ba\'ia, Matias, and Santos \cite{bms},
confirms the roles of $\left( trM\right) ^{+}$, $\left( trM\right) ^{-}$,
and $\left\vert trM\right\vert $ established by Owen and Paroni. \ In doing
so, we give an alternative, shorter proof of Owen and Paroni's results, and
we establish additional explicit formulas for other measures of
disarrangements. \

In order to motivate our study and to provide necessary background,
we discuss briefly in the following subsections of this introduction some
concepts and results from the multiscale geometry of structured
deformations. \ (Readers familiar with this material may wish to skip to the
last subsection of the introduction where our main results are summarized.)
\ 

\subsection{Structured deformations and disarrangement densities in the setting of Del Piero and Owen}

The need in continuum mechanics to include the effects of multiscale
geometrical changes led Del Piero and Owen \cite{del piero owen} to a notion
of structured deformations as triples $(\kappa ,g,G)$, where

\begin{itemize}
\item the injective, piecewise continuously differentiable field $g$ maps
the points of a continuous body into physical space and describes
macroscopic changes in the geometry of the body,

\item the piecewise continous tensor field $G$ maps the body into the space
of linear mappings on the translation space of physical space and satisfies
the "accommodation inequality"%
\begin{equation}
0<C<\det G(x)\leq \det \nabla g(x)\text{ \ at each point }x
\label{accommodation}
\end{equation}%
where $\nabla $ denotes the classical gradient operator, and

\item $\kappa $ is a surface-like subset of the body that describes
preexisting, unopened macroscopic cracks.
\end{itemize}

\  A geometrical interpretation of the field $G$ is provided by the
Approximation Theorem \cite{del piero owen}: for each structured deformation 
$(\kappa ,g,G)$ there exists a sequence of injective, piecewise smooth
deformations $f_{n}$ and a sequence of surface-like subsets $\kappa _{n}$ of
the body such that%
\begin{equation}
g=\lim_{n\to \infty }f_{n}  \label{fn converges to g}
\end{equation}%
\begin{equation}
G=\lim_{n\to \infty }\nabla f_{n}  \label{grad fn converges to G}
\end{equation}%
and%
\begin{equation}
\kappa =\bigcup _{n=1}^{\infty }\bigcap _{p=n}^{\infty }\kappa _{p}
\label{kappa n converges to kappa}
\end{equation}%
The limits in \eqref{fn converges to g} and \eqref{grad fn converges to G}
are taken in the sense of $L^{\infty }$convergence. A sequence $n\longmapsto
f_{n}$ of piecewise smooth, injective functions satisfying \eqref{fn converges to g} and \eqref{grad fn converges to G} is called a \textit{%
determining sequence} for the pair $(g,G)$, and each term $f_{n}$ is
interpreted as describing the body divided into tiny pieces that may
individually undergo smooth geometrical changes and that also may undergo 
\textit{disarrangements}, i.e., may separate or slide relative to each
other. \ In this context, we write $f_{n}\rightsquigarrow (g,G)$. \ From
\eqref{grad fn converges to G} we see that $G$ captures the effects at the
macrolevel of smooth geometrical changes at submacroscopic levels, and we
call $G$ the \textit{deformation without disarrangements}. \ 

Del Piero and Owen \cite{integral gradient} proved that, for every
structured deformation $(\kappa ,g,G)$, for every determining sequence $%
n\longmapsto f_{n}$ for $(g,G)$, and for every point $x$ where $g$ is
differentiable and where $G$ is continuous, there holds%
\begin{equation}
\lim_{r\to 0}\lim_{n\to \infty }\frac{%
\displaystyle \int_{J(f_{n})\cap {B}_{r}(x)}[f_{n}](y)\otimes \nu (y)\,d\mathcal{H}%
^{N-1}(y)}{\left\vert B_{r}(x)\right\vert }=\nabla g(x)-G(x)\text{.}
\label{M as a limit of jumps}
\end{equation}%
Here, $\mathcal{H}^{N-1}$ denotes the $(N{-}1)$-dimensional Hausdorff measure
on $\mathbb{R}^{N}$, $B_{r}(x)$ denotes the open ball centered at $x$ of radius $r
$, $\left\vert B_{r}(x)\right\vert $ denotes its volume (i.e., its $N$-dimensional Lebesgue measure), $J(f_{n})$ denotes the jump set of $f_{n}$,
i.e., points where $f_{n\text{ }}$ can suffer jump-discontinuities, and $%
[f_{n}](y)\otimes \nu (y)$ is the tensor product of the jump $[f_{n}]$ of $%
f_{n}$ with the normal $\nu \,$\ to the jump set. \ This result permits us
to call the tensor
\begin{equation}\label{definition of M}
M(x):=\nabla g(x)-G(x)  
\end{equation}
the \textit{deformation due to disarrangements}, because it captures, in the
limit as $n$ tends to infinity, the volume density of separations and slips
between pieces of the body as described by the approximating deformations $%
f_{n}$. \ We may then regard the tensor field $M$ as a \textit{tensorial
disarrangement density} that, for every determining sequence $n\longmapsto
f_{n}$ for $(g,G)$, reflects the limits of interfacial discontinuities of
the approximating deformations $f_{n}$. Moreover, \eqref{fn converges to g}
and \eqref{grad fn converges to G} along with the definition of $M$ \eqref{definition of M} yield the alternative formula for the disarrangement
density:%
\begin{equation}
M=\nabla (\lim_{n\to \infty }f_{n})-\lim_{n\to\infty }\nabla f_{n}.  \label{gradient of limit minus limit of gradient}
\end{equation}%
Consequently, $M$ measures quantitatively the lack of commutativity of the
classical gradient $\nabla $ and the limit operator $\lim_{n\longrightarrow
\infty }$ for $L^{\infty }$-convergence. \ 

The trivial algebraic relation%
\begin{equation}
\nabla g=G+M  \label{additive decomposition of grad g}
\end{equation}%
together with the identification relations \eqref{grad fn converges to G}
and \eqref{M as a limit of jumps} shows that the macroscopic deformation
gradient $\nabla g$ has an additive decomposition into its part $G$ without
disarrangements and its part $M$ due to disarrangements. \ Because $G$ has
invertible values, \eqref{additive decomposition of grad g} leads
immediately to two multiplicative decompositions for $\nabla g$:%
\begin{equation}
\nabla g=G(I+G^{-1}M)=(I+MG^{-1})G.
\label{multiplicative decompositions for grad g}
\end{equation}

The disarrangement density $M$ and the deformation without disarrangements $G
$ have an additional property significant in the description of defects and
dislocations in a continuous body in three dimensions. We consider a smooth
surface $\mathcal{S}$ with smooth bounding closed curve $\gamma$, both contained
in a region in the body where $g$ and $G$ are smooth. The relation \eqref{additive decomposition of grad g} and the smoothness of $g~$imply%
\begin{equation*}
0=\oint_{\gamma}\nabla g(x)dx=\oint_{\gamma}G(x)dx+\oint_{\gamma}M\left( x\right) dx.
\end{equation*}%
The vector $\oint_{\gamma}M\left( x\right) dx$ measures the displacement due to
disarrangements along $\gamma$ and may be called the Burgers vector \cite{del
piero owen} for $\gamma$ arising from the given structured deformation. \
Application of Stokes' Theorem to $\oint_{\gamma}G(x)dx$ and $\oint_{\gamma}M\left(
x\right) dx$ and use of the previous relation yields the formulas for the
Burgers vector:%
\begin{equation}
\oint_{\gamma}M\left( x\right) dx=\int_{\mathcal{S}}\mathrm{curl}M(x)\nu
(x)dA_{x}=-\int_{\mathcal{S}}\mathrm{curl}G(x)\nu (x)dA_{x}.
\label{formula for burgers vector}
\end{equation}%
The second-order tensor field $\mathrm{curl}M=-\mathrm{curl}G$ thus determines
the Burgers vector associated with $\gamma$ for every closed curve and
corresponds to familiar measures of dislocation density \cite{Kroner,Nye}. 
In this manner, the disarrangement density tensor $M$ determines both the Burgers vector and
the dislocation density tensor, both basic tools in modelling the effects of
submacroscopic defects on the response of solids. \ \ \ \ 

The tensorial relations \eqref{definition of M} and \eqref{M as a limit of jumps} yield upon application of the trace operator the scalar relation%
\begin{equation}
\lim_{r\to 0}\lim_{n\to \infty }\frac{\displaystyle
\int_{J(f_{n})\cap {B}_{r}(x)}[f_{n}](y)\cdot \nu (y)\,d\mathcal{H}%
^{N-1}(y)}{\left\vert B_{r}(x)\right\vert }=trM(x)
\label{trM as a limit of normal jumps}
\end{equation}%
in which $[f_{n}](y)\cdot \nu (y)$ is the scalar product of the jump and of
the normal at $y$. \ The formula \eqref{trM as a limit of normal jumps}
tells us that $trM$ is a \textit{scalar (bulk) disarrangement density} that
captures the components of the jumps of $f_{n}$ that are normal to the jump
set. Moreover, this scalar disarrangement density at $x$, $trM(x),$ allows
for cancellation of positive and negative contributions of $[f_{n}](y)\cdot
\nu (y)$ at points $y$ near $x$ to the integral on the left-hand side of \eqref{trM as a limit of normal jumps}. Thus, $trM(x)$ does not distinguish
between jumps with $[f_{n}](y)\cdot \nu (y)>0$ that pull apart small pieces
of the body near $x$ and jumps with $[f_{n}](y)\cdot \nu (y)<0$ that cause
small pieces near $x$ to switch places. Because the approximating $f_{n}$
are injective, the possibility for the case $[f_{n}](y)\cdot \nu (y)<0$ that 
$f_{n}$ can cause adjacent small pieces of the body to interpenetrate is
ruled out.

Owen and Paroni \cite{OwenParoni} refined the scalar disarrangement density $%
trM$ by replacing $[f_{n}](y)\cdot \nu (y)$ by its positive part throughout
the jump set of $f_{n}$, or by its negative part throughout the jump set:%
\begin{equation}
([f_{n}](y)\cdot \nu (y))^{+}=\tfrac{1}{2}(\left\vert [f_{n}](y)\cdot \nu
(y)\right\vert +[f_{n}](y)\cdot \nu (y))  \label{positive part}
\end{equation}%
\begin{equation}
([f_{n}](y)\cdot \nu (y))^{-}=\tfrac{1}{2}(\left\vert [f_{n}](y)\cdot \nu
(y)\right\vert -[f_{n}](y)\cdot \nu (y)).  \label{negative part}
\end{equation}%
The field $([f_{n}]\cdot \nu \,)^{+}$ on the jump set is a scalar
(interfacial) disarrangement density that measures separations of pieces of
the body caused by $f_{n}$, while the field $([f_{n}]\cdot \nu \,)^{-}$ is a
scalar (interfacial) disarrangement density that measures the switching of
pieces of the body caused by $f_{n}$. Since we have%
\begin{equation*}
\left\vert \lbrack f_{n}](y)\cdot \nu (y)\right\vert =([f_{n}](y)\cdot \nu
(y))^{+}+([f_{n}](y)\cdot \nu (y))^{-}\text{,}
\end{equation*}%
the field $\left\vert [f_{n}](y)\cdot \nu (y)\right\vert $ is a scalar
disarrangement density that measures both separations and switches. We fix a
part $\mathcal{P}$ of the body, we integrate \eqref{positive part} or \eqref{negative part} over $J(f_{n})\cap \mathcal{P}$\ and use the formula \eqref{trM as a limit of normal jumps} to obtain the relations 
\begin{equation}\label{lim infs of positive and negative parts}
\begin{split}
\liminf_{n\to \infty} &\int_{J(f_{n})\cap \mathcal{P}}  ([f_{n}](y)\cdot \nu (y))^{\pm }d\mathcal{H}^{N-1}(y) \\
=& \frac{1}{2}\liminf_{n\to \infty }\int_{J(f_{n})\cap \mathcal{%
P}}\left\vert [f_{n}](y)\cdot \nu (y)\right\vert d\mathcal{H}^{N-1}(y) \\
& \pm \frac{1}{2}\liminf_{n\to \infty }\int_{J(f_{n})\cap 
\mathcal{P}}[f_{n}](y)\cdot \nu (y)d\mathcal{H}^{N-1}(y)   \\
=&\frac{1}{2}\liminf_{n\to \infty }\int_{J(f_{n})\cap \mathcal{%
P}}\left\vert [f_{n}](y)\cdot \nu (y)\right\vert d\mathcal{H}^{N-1}(y)
\pm \frac{1}{2}\int_{\mathcal{P}}trM(x)d\mathcal{L}^{N}(x).
\end{split}
\end{equation}%
Consequently, the limiting behavior of the integral of $([f_{n}](y)\cdot \nu
(y))^{\pm }$ in \eqref{lim infs of positive and negative parts} as $n$ tends
to $\infty $ is determined by the behavior of the integral of $\left\vert
[f_{n}](y)\cdot \nu (y)\right\vert $, and we restrict our attention to the
latter. \ We expect that%
\begin{equation*}
\liminf_{n\to \infty }\int_{J(f_{n})\cap \mathcal{P}}\left\vert
[f_{n}](y)\cdot \nu (y)\right\vert d\mathcal{H}^{N-1}(y),
\end{equation*}%
unlike%
\begin{equation*}
\liminf_{n\to \infty }\int_{J(f_{n})\cap \mathcal{P}%
}[f_{n}](y)\cdot \nu (y)d\mathcal{H}^{N-1}(y),
\end{equation*}%
will depend upon the choice of determining sequence for $(g,G)$. \
Therefore, we are led to consider the most economical manner in which
separations and switches can arise among the determining sequences for $(g,G)
$:%
\begin{equation}
\mathcal{V}^{\left\vert \cdot \right\vert }(g,G;\mathcal{P}):=\inf \left\{
\liminf_{n\to \infty }\int_{J(f_{n})\cap \mathcal{P}}\left\vert
[f_{n}](y)\cdot \nu (y)\right\vert {\small d}\mathcal{H}^{N-1}{\small %
(y)\,\,:\,}f_{n}\rightsquigarrow (g,G)\right\} .
\label{minimal disarrangements in P}
\end{equation}%
The number $\mathcal{V}^{\left\vert \cdot \right\vert }(g,G;\mathcal{P})$ so
defined has the dimension of volume, and we call $\mathcal{V}^{\left\vert
\cdot \right\vert }(g,G;\mathcal{P})$ the \textit{(minimal) volume swept out
by disarrangements in }$\mathcal{P}$\textit{\ for }$(g,G)$. If we replace $%
\left\vert \cdot \right\vert $ everywhere in \eqref{minimal disarrangements in P} by $"+"$ or everywhere by $"-"$, then we call the number\ $\mathcal{V}%
^{+}(g,G;\mathcal{P})$ the \textit{(minimal) volume swept out by separations
in }$\mathcal{P}$\textit{\ for }$(g,G)$, and the number $\mathcal{V}^{-}(g,G;%
\mathcal{P})$ the \textit{(minimal) volume swept out by switches in} $%
\mathcal{P}$\textit{\ for }$(g,G)$. \ The formulas \eqref{lim infs of positive and negative parts} imply the simple formulas%
\begin{equation}
\mathcal{V}^{\pm }(g,G;\mathcal{P})=\tfrac{1}{2}\mathcal{V}^{\left\vert
\cdot \right\vert }(g,G;\mathcal{P})\pm \tfrac{1}{2}\int_{\mathcal{P}}trM(x)d%
\mathcal{L}^{N}(x)  \label{formulas for v plus or minus}
\end{equation}%
and, in view of the form of the second term on the right-hand side, raise the
following basic question: \ does the volume swept out by disarrangements $%
\mathcal{V}^{\left\vert \cdot \right\vert }(g,G;\mathcal{P})$ as defined in \eqref{minimal disarrangements in P} have an associated disarrangement
density which, when integrated over $\mathcal{P}$, recovers $\mathcal{V}%
^{\left\vert \cdot \right\vert }(g,G;\mathcal{P})$. If so, what specific
information can be obtained about the dependence of the integrand upon the
structured deformation $(g,G)$?

While the setting for structured deformations described in this subsection is
quite suitable for formulating refined field equations in continuum
mechanics \cite{DeseriOwen} that reflect the influence of submacroscopic geometrical changes
in a body, this setting has not provided answers to questions such as the
ones just raised. \ Part of the difficulty with the setting provided in \cite%
{del piero owen} lies in the choice of smoothness placed on $g$ and its
approximates $f_{n}$, while another part lies in the requirement that $g$
and $f_{n}$ be injective. \ An alternative setting provided by Choksi and
Fonseca \ \cite{choksifonseca} \ was proposed for dealing with such
questions and is described briefly in the next subsection.

\subsection{Structured deformations and disarrangement densities in the setting of Choksi and Fonseca}

We describe here a few essential elements of the treatment of structured
deformations by Choksi and Fonseca \ \cite{choksifonseca}. The articles \cite%
{cdfo}, \cite{bms}, \cite{bms survey}, and \cite{mirek} also provide
summaries of that treatment, and \cite{bms}, \cite{bms survey}, and \cite%
{mirek} provide alternative settings for structured deformations. The
summary in \cite{cdfo} is intended for those interested in immediate
applications in continuum mechanics, while \cite{bms} sets the stage for
applications of structured deformations to thin bodies \cite{MatiasSantos}.
\ The article \cite{mirek} reexamines the results of \cite{choksifonseca} in
a broader setting while providing refinements of counterparts of the
Approximation Theorem and the identification relation \eqref{M as a limit of jumps}. 

According to Choksi and Fonseca, a structured deformation is a pair $(g,G)$
in which $g:\Omega \longrightarrow \mathbb{R}^{N}$, with $\Omega $ an open
subset of the space $\mathbb{R}^{N}$ of $N$-tuples of real numbers, and $G$ $%
:\Omega \longrightarrow \mathbb{R}^{N\times N}$, with $\mathbb{R}^{N\times N}$%
the space of $N{\times}N$ matrices with real entries. \ The mapping $G$ is
assumed to be integrable on $\Omega $, $G\in L^{1}(\Omega;\mathbb{R}%
^{N\times N})$, and $g$ is assumed to be in the space $SBV(\Omega;\mathbb{R}%
^{N})$, i.e., $g$ is a function of bounded variation with the additional
property that its distributional derivative $Dg$, as a bounded measure, has
zero Cantor part: 
\begin{equation}
Dg=\nabla g\,\mathcal{L}^{N}+[g]\otimes \nu \mathcal{H}^{N-1}.
\label{distributional derivative of g}
\end{equation}%
Here the integrable mapping $\nabla g$ is the density of the absolutely
continuous part $\nabla g\,\mathcal{L}^{N}$of $Dg$ with respect to $N$%
-dimensional Lebesgue measure $\mathcal{L}^{N}$, and $[g]\otimes \nu $ is
the density of the singular part $[g]\otimes \nu \mathcal{H}^{N-1}$of $Dg$
with respect to $(N{-}1)$-dimensional Hausdorff measure $\mathcal{H}^{N-1}$. \
The singular part is concentrated on $J(g)$, the jump set of $g$, and, as
usual, $[g]$ denotes the jump in $g$ and $\nu $ denotes the normal to the
jump set $J(g)$. It is important to note that $\nabla g$ in the present
setting is no longer the classical gradient of a smooth field and,
consequently, need not be curl-free. \ Nevertheless, $\nabla g$ satisfies an
integral version of the property of approximation by linear mappings that
defines the classical gradient of smooth fields .

Choksi and Fonseca \cite{choksifonseca} prove a version of the Approximation
Theorem with approximating deformations $f_{n}$ also in $SBV(\Omega;\mathbb{%
R}^{N})$ and with \eqref{fn converges to g} and \eqref{grad fn converges to G} replaced respectively by%
\begin{equation}
f_{n}\to g\text{ in }L^{1}(\Omega;\mathbb{R}^{N})\text{ \ }
\label{L1 convergence of fn to g}
\end{equation}%
and%
\begin{equation}
\nabla f_{n}\rightharpoonup G\text{ \ weakly in the sense of measures.}
\label{weak convergence of grad fn to G}
\end{equation}%
We note that no restriction in the form of the accommodation inequality \eqref{accommodation} or in the form of a requirement of injectivity of $g$ or $%
f_{n}$ is imposed in the present context. \ We again use the term
\emph{determining sequence} to describe a sequence $n\longmapsto f_{n}$ satisfying \eqref{L1 convergence of fn to g} and \eqref{weak convergence of grad fn to G} for a given structured deformation $(g,G)$, and we again write $%
f_{n}\rightsquigarrow (g,G)$ when \eqref{L1 convergence of fn to g} and \eqref{weak convergence of grad fn to G} both hold. \ The properties of
distributional derivatives along with relations \eqref{distributional derivative of g}, \eqref{L1 convergence of fn to g}, and \eqref{weak convergence of grad fn to G} justify the calculation%
\begin{eqnarray*}
\nabla g\,\mathcal{L}^{N}+[g]\otimes \nu \,\mathcal{H}^{N-1}
&=&D\lim_{n\rightarrow \infty }f_{n} \\
&=&\lim_{n\rightarrow \infty }Df_{n} \\
&=&\lim_{n\rightarrow \infty }(\nabla f_{n}\,\mathcal{L}^{N}+[f_{n}]\otimes
\nu \,\mathcal{H}^{N-1}) \\
&=&G\,\mathcal{L}^{N}+\lim_{n\rightarrow \infty }([f_{n}]\otimes \nu \,%
\mathcal{H}^{N-1})
\end{eqnarray*}%
where the convergence indicated in the last three lines is weak convergence
in the sense of measures. We conclude that the singular parts $%
[f_{n}]\otimes \nu \,\mathcal{H}^{N-1}$ of the approximating deformations $%
f_{n}$ converge in the same sense and that their limit satisfies%
\begin{equation}
\lim_{n\rightarrow \infty }([f_{n}]\otimes \nu \,\mathcal{H}^{N-1})=(\nabla
g\,-G)\mathcal{L}^{N}+[g]\otimes \nu \,\mathcal{H}^{N-1}.
\label{limit of singular parts}
\end{equation}%
In particular, the restriction of the limiting measure $\lim_{n\rightarrow
\infty }([f_{n}]\otimes \nu \,\mathcal{H}^{N-1})$ to the complement of the
jump set $J(g)$ agrees with the corresponding restriction of \ $(\nabla
g\,-G)\mathcal{L}^{N}=M\,\mathcal{L}^{N}$. \ Consequently, the tensor
field $M=$\ $\nabla g\,-G$ retains in this broader setting for structured
deformations its identity as a tensor density of disarrangements for $(g,G)$%
. \ The\ formula \eqref{limit of singular parts} shows that when $M=\nabla
g-G\neq 0$, while all of the measures $[f_{n}]\otimes \nu \,\mathcal{H}%
^{N-1} $ are supported on sets $J(f_{n})$ of $\mathcal{L}^{N}$-measure zero
and so have $\mathcal{L}^{N}$-parts zero, the limit measure $%
\lim_{n\rightarrow \infty }([f_{n}]\otimes \nu \,\mathcal{H}^{N-1})$ has $%
\mathcal{L}^{N}$-part $M\mathcal{L}^{N}$non-zero. \ This observation points
to the fact that the jump sets $J(f_{n})$ can diffuse in the limit
throughout the domain $\Omega $ so that the limiting measure $%
\lim_{n\rightarrow \infty }([f_{n}]\otimes \nu \,\mathcal{H}^{N-1})$ is
supported in part on sets of positive $\mathcal{L}^{N}$-measure. \ This
provides a counterpart in the $SBV$-setting to the relation \eqref{M as a limit of jumps} in which the limit of jumps on the left-hand side delivers
the $\mathcal{L}^{N}$-density $M$. \ (See \cite{mirek} for a detailed
derivation of a counterpart of \eqref{M as a limit of jumps} in a somewhat
broader setting for structured deformations than $SBV$.)

We note briefly that the scalar density of disarrangements $trM=$ $tr(\nabla
g\,-G)$ that counts only normal components of jumps and that emerged in the
previous setting also appears in the present setting when one takes the
trace of every member of \eqref{limit of singular parts}: if $%
f_{n}\rightsquigarrow (g,G)$, then 
\begin{equation}
\lim_{n\rightarrow \infty }([f_{n}]\cdot \nu \,\mathcal{H}^{N-1})=tr(\nabla
g\,-G)\mathcal{L}^{N}+[g]\cdot \nu \,\mathcal{H}^{N-1}.
\label{trace M as a scalar density}
\end{equation}%
However, as was the case in the setting of Del Piero and Owen, replacement
of $[f_{n}]\cdot \nu $ by $([f_{n}]\cdot \nu )^{\pm }$ or by $\left\vert
[f_{n}]\cdot \nu \right\vert $ need not yield a limit of the corresponding
measures and, if a limit exists, the limit may depend upon the choice of
determining sequence $n\longmapsto f_{n}$. \ The setting of Choksi and
Fonseca was formulated as a means of resolving these difficulties, and we
summarize some aspects of that resolution in the next subsection.

\subsection{Relaxation of energies for structured deformations}
\label{sss}

Optimal functions arising from structured deformations such as the one \eqref{minimal disarrangements in P}
\begin{equation*}
\mathcal{V}^{\left\vert \cdot \right\vert }(g,G;\mathcal{P})=\inf \bigg\{
\liminf_{n\longrightarrow \infty }\int_{J(f_{n})\cap \mathcal{P}}\left\vert
[f_{n}](y)\cdot \nu (y)\right\vert {\small d}\mathcal{H}^{N-1}{\small %
(y)\,\,:\,}f_{n}\rightsquigarrow (g,G)\bigg\} 
\end{equation*}%
introduced in Section 1.1 can be analyzed using the results of Choksi and
Fonseca \cite{choksifonseca} on "relaxation of energies" for structured
deformations. \ In that approach, the integral $\int_{J(f_{n})\cap \mathcal{P%
}}\left\vert [f_{n}](y)\cdot \nu (y)\right\vert {\small d}\mathcal{H}^{N-1}%
{\small (y)}$ is replaced by an initial energy functional 
\begin{equation}
E(f_{n})=\int_{\Omega }W(\nabla f_{n}(y)){\small d}\mathcal{L}^{N}{\small (y)%
}+\int_{J(f_{n})\cap \Omega }\psi ([f_{n}](y),\nu (y)){\small d}\mathcal{H}%
^{N-1}{\small (y)}  \label{initial energy functional}
\end{equation}%
defined for $f_{n}\in SBV(\Omega;\mathbb{R}^{N})$. \ By imposing conditions
on the initial bulk energy density $W$ and on the initial interfacial energy
density $\psi $, the goal is to obtain for the relaxed energy $I(g,G)$
defined by 
\begin{equation}\label{definition of relaxed energy}
\begin{split}
I(g,G):=\inf\bigg\{ \liminf_{n\to \infty}\bigg(& \int_{\Omega}%
W(\nabla f_{n}(y))d\mathcal{L}^{N}(y)+ \\
&  \int_{J(f_{n})\cap \Omega } \psi ([f_{n}](y),\nu (y){\small )d}\mathcal{H}^{N-1}{\small (y)}%
\bigg): f_{n}\rightsquigarrow (g,G) \bigg\}
\end{split}
\end{equation}%
a representation of the form 
\begin{equation}
I(g,G)=\int_{\Omega }H(\nabla g(y),G(y)){\small d}\mathcal{L}^{N}{\small (y)}%
+\int_{J(g)\cap \Omega }h([g](y),\nu (y)){\small d}\mathcal{H}^{N-1}{\small %
(y)}  \label{representation for relaxed energy}
\end{equation}%
and to deduce properties of the relaxed bulk energy density $H$ and the
relaxed interfacial energy density $h$. Because our present interest lies in
the case of disarrangement densities, and not on the full energetics of
structured deformations, we shall restrict our attention to the case $W=0$, and
we record the\ following adaptation for the case $W=0$ of results from \cite%
{choksifonseca} (see \cite[Theorem 3]{OwenParoni} for further comments and
other adaptations).

\begin{theorem}\label{th1}
Let $S^{N-1}=\{\nu \in \mathbb{R}^{N}:|\nu |=1\}.$ Let $\Omega $ be a
bounded open subset of $\mathbb{R}^{N}$ and $\psi :\mathbb{R}^{N}\times
S^{N-1}\rightarrow \lbrack 0,+\infty )$ be such that

\noindent (H1) there exists a constant $C>0$ such that 
\begin{equation}
0\leq \psi (\xi ,\nu )\leq C|\xi |  \label{growth of psi}
\end{equation}%
for all $(\xi ,\nu )\in \mathbb{R}^{N}\times S^{N-1}$,

\noindent (H2) $\psi (\cdot ,\nu )$ is positively homogeneous of degree 1: 
\begin{equation}
\psi (t\,\xi ,\nu )=t\,\psi (\xi ,\nu )  \label{homogeneity of psi}
\end{equation}%
for all $t>0$ and $(\xi ,\nu )\in \mathbb{R}^{N}\times S^{N-1}$,

\noindent (H3) $\psi(\cdot,\nu)$ is subadditive, i.e., for all $\xi _{1},\xi _{2}\in 
\mathbb{R}^{N}$ and $\nu \in S^{N-1}$, 
\begin{equation}
\psi (\xi _{1}+\xi _{2},\nu )\leq \psi (\xi _{1},\nu )+\psi (\xi _{2},\nu ).
\label{subadditivity of psi}
\end{equation}%
Then, for any $p>1$, if we define 
\begin{equation*}
\begin{split}
I(g,G):=\inf\bigg\{ &\liminf_{n\rightarrow \infty}\int_{J(u_{n})\cap \Omega }\psi ([u_{n}],\nu )\,d\mathcal{H}^{N-1}:u_{n}\in
SBV(\Omega;\mathbb{R}^{N}), \\
&\; u_{n}\rightarrow g\mbox{ in }L^{1}(\Omega;\mathbb{R}^{N}),\nabla u_{n}\overset{\ast }{\rightharpoonup }G, \\
&\; \sup_{n}\left( |\nabla u_{n}|_{L^{p}(\Omega;\mathbb{R}^{N\times N})}+|Du_{n}|(\Omega )\right) <+\infty \bigg\},
\end{split}
\end{equation*}%
we have 
\begin{equation*}
I(g,G)=\int_{\Omega }H(\nabla g(x),G(x))\,d\mathcal{L}^{N}+\int_{J(g)\cap
\Omega }h([g](x),\nu (x))\,d\mathcal{H}^{N-1}(x),
\end{equation*}%
where 
\begin{equation}\label{relaxed bulk density}
\begin{split}
H(A,B):= \inf\bigg\{ & \int_{J(u)\cap Q}\psi ([u],\nu )\,d\mathcal{H}%
^{N-1}:u\in SBV(Q;\mathbb{R}^{N}),  \\
&\; u|_{\partial Q}=Ax,|\nabla u|\in L^{p}(Q),\int_{Q}\nabla u\,d%
\mathcal{L}^{N}=B\bigg\},  
\end{split}
\end{equation}%
and 
\begin{equation}\label{relaxed interfacial density}
\begin{split}
h(\xi ,\eta ):=\inf\bigg\{ & \int_{J(u)\cap Q_{\eta }}\psi ([u],\nu )\,d%
\mathcal{H}^{N-1} :u\in SBV(Q_{\eta }; \mathbb{R}^{N}),  \\
&\; u|_{\partial Q_{\eta}}=u_{\xi ,\eta }, \nabla u=0\;\text{a.e.}\bigg\},  
\end{split}
\end{equation}%
with 
\begin{equation}
u_{\xi ,\eta }(x):=\left\{ 
\begin{array}{cl}
0 & \mbox{ if }~-\frac{1}{2}\leq x\cdot \eta <0, \\ 
\xi & \mbox{ if }~0\leq x\cdot \eta <\frac{1}{2}.%
\end{array}%
\right.  \label{uxinu}
\end{equation}%
Here, $Q=(-1/2,1/2)^{N}$ and $Q_{\eta }$ denotes the unit cube
centered at the origin and with two faces normal to $\eta $.
\end{theorem}
In the right-hand side of \eqref{relaxed interfacial density} we have corrected an inconsequential misprint that is present in
the corresponding formula in Theorem 3 of \cite{OwenParoni}.

Another approach to relaxation of energies for structured deformations in the full $BV$ setting is provided in \cite{bms}. 
A structured deformation in \cite{bms} is a pair $(g,G)\in BV^2(\Omega;\mathbb{R}^{N})\times BV(\Omega;\mathbb{R}^{N \times N})$, where $ BV^2(\Omega;\mathbb{R}^{N}) :=  \{ u \in BV(\Omega;\mathbb{R}^{N}) : \nabla u \in BV(\Omega;\mathbb{R}^{N \times N})\}$.
The counterpart of the Approximation Theorem in this context asserts that there exists a sequence $f_n \in BV^2(\Omega;\mathbb{R}^{N})$ such that both $f_n \to g$ and $\nabla f_n \to G$ in the $L^1$-norm.
In this case we write $f_{n}\rightsquigarrow (g,G)$.

The energy functional considered in \cite{bms}, under assumptions on the initial bulk and surface energy densities similar to the ones in \cite{choksifonseca}, reads
\begin{equation}\label{bms-funct}
\begin{split}
E(f_n)=&\int_\Omega W(\nabla f_n(y),\nabla^2 f_n(y))\,d\mathcal{L}^N{y}  +\int_{J(f_n)} \psi([f_n](y),\nu(y))\,d\mathcal{H}^{N-1}(y) \\
& +\int_{J(\nabla f_n)} \psi_1([\nabla f_n](y),\nu(y))\, d\mathcal{H}^{N-1}(y),
\end{split}
\end{equation}
\noindent and the relaxed energy $I(g,G)$ is 
defined by 
\begin{equation}\label{bms-rel-ene}
I(g,G):=\inf\Big\{\liminf_{n\to\infty} E(f_n): f_{n}\rightsquigarrow (g,G)\Big\}.
\end{equation}

A crucial result in \cite{bms} is that \eqref{bms-rel-ene} can be divided into two first-order relaxed energies, namely, $I(g,G) = I_1(g, G) + I_2(G)$, where the term $I_1(g, G)$ captures the structured deformation, whereas $I_2(G)$ only depends on the deformation without disarrangements $G$.
In the relevant case for the present paper, i.e., $W=\psi_1=0$, the results in \cite{bms} give $I_2=0$ and 
\begin{equation}\label{I_1}
I_1(g,G):= \inf\bigg\{ \liminf_{n\to \infty}  \int_{J(f_{n})\cap \Omega } \psi ([f_{n}](y),\nu (y){\small )d}\mathcal{H}^{N-1}{\small (y)}: f_{n}\rightsquigarrow (g,G)\bigg\}. 
\end{equation}

Defining $SBV^2(\Omega;\mathbb{R}^{N}) :=  \{ u \in SBV(\Omega;\mathbb{R}^{N}) : \nabla u \in SBV(\Omega;\mathbb{R}^{N \times N})\}$, the following representation theorem holds
\begin{theorem}[see {\cite[Theorem 3.2]{bms}}]\label{bms-thm}
For every $(g,G)\in SBV^2(\Omega;\mathbb{R}^{N})\times SBV(\Omega;\mathbb{R}^{N \times N})$, 
given $\psi$ under the same hypotheses (H1)-(H3) of Theorem \ref{th1}, we have that
\begin{equation*}
I(g,G)=\int_{\Omega }H(G(x) - \nabla g(x))\,d\mathcal{L}^{N}+\int_{J(g)\cap
\Omega }h([g](x),\nu (x))\,d\mathcal{H}^{N-1}(x),
\end{equation*}
where, given $A\in \R{N\times N}$, $\xi\in\R N$, and $\eta\in S^{N-1}$, 
\begin{equation}\label{relaxed bulk density in BMS}
\begin{split}
H(A):= \inf\bigg\{ & \int_{J(u)\cap Q}\psi ([u],\nu )\,d\mathcal{H}%
^{N-1}:u\in SBV^2(Q;\mathbb{R}^{N}),  \\
&\; u|_{\partial Q}= 0, \nabla u = A  \; \text{a.e. in } Q \bigg\}
\end{split}
\end{equation}
and 
\begin{equation}\label{relaxed interfacial density in BMS}
\begin{split}
h(\xi ,\eta ):=\inf\bigg\{ & \int_{J(u)\cap Q_{\eta }}\psi ([u],\nu )\,d\mathcal{H}^{N-1} :u\in SBV^2(Q_{\eta }; \mathbb{R}^{N}),  \\
&\; u|_{\partial Q_{\eta}}=u_{\xi ,\eta }, \nabla u=0\;\text{a.e. in } Q \bigg\},  
\end{split}
\end{equation}%
with $u_{\xi,\eta}$ defined as in \eqref{uxinu}.
\end{theorem}
\begin{remark}\label{bms-remark}
It is worth noticing that the minimum problems defining \eqref{relaxed bulk density in BMS} and \eqref{relaxed interfacial density in BMS} are formally performed in $SBV^2(\Omega;\mathbb{R}^{N})$, but the result is the same if $SBV^2$ is replaced in these relations by $SBV$, due to the requirement that $\nabla u$ be constant.
\end{remark}

\subsection{Explicit formulas for relaxed disarrangement densities}

Owen and Paroni \cite{OwenParoni} applied Theorem \ref{th1} to the specific
disarrangement densities $\left\vert [f_{n}](y)\cdot \nu (y)\right\vert $
and $([f_{n}](y)\cdot \nu (y))^{\pm }$ introduced in Section 1.1 and
obtained for each of these densities an explicit formula for the
corresponding relaxed disarrangement densities $H$ in \eqref{relaxed bulk density} and $h$ in \eqref{relaxed interfacial density}. \ Among their results
(\cite{OwenParoni},Theorem 4) are the following (obtained by setting $L(x)=I$
in their Theorem 4):

\begin{theorem}\label{th2}
The initial disarrangement densities%
\begin{equation}
\psi ^{\left\vert \cdot \right\vert }(\xi ,\nu ):=\left\vert \xi \cdot \nu
\right\vert  \label{psi absolute}
\end{equation}%
\begin{equation}
\psi ^{\pm }(\xi ,\nu ):=(\xi \cdot \nu )^{\pm }  \label{psi plus or minus}
\end{equation}%
satisfy the hypotheses (H1)-(H3) in Theorem \ref{th1} and have relaxed
disarrangement densities given by%
\begin{equation}
H^{\left\vert \cdot \right\vert }(A,B)=\left\vert tr(A-B)\right\vert \text{,
\ \ }h^{\left\vert \cdot \right\vert }(\xi ,\nu )=\left\vert \xi \cdot \nu
\right\vert =\psi ^{\left\vert \cdot \right\vert }(\xi ,\nu ),
\label{absolute value relaxed densities}
\end{equation}%
and%
\begin{equation}
H^{\pm }(A,B)=(tr(A-B))^{\pm }\text{, \ \ }h^{\pm }(\xi ,\nu )=(\xi \cdot
\nu )^{\pm }=\psi ^{\pm }(\xi ,\nu ).
\label{plus or minus relaxed densities}
\end{equation}
\end{theorem}

Specifically, when the minimal volume swept out by disarrangements $\mathcal{%
V}^{\left\vert \cdot \right\vert }(g,G;\mathcal{P})$ is defined in the
Choksi-Fonseca setting by \eqref{minimal disarrangements in P}, then \eqref{absolute value relaxed densities} yields the explicit formula%
\begin{equation}\label{formula for volume swept out}
{\small \mathcal{V}}^{\left\vert \cdot \right\vert }{\small (g,G;\mathcal{P})} =\int_{%
\mathcal{P}}\left\vert tr(\nabla g(x)-G(x))\right\vert {\small d}\mathcal{L}%
^{N}{\small (x)+}\int_{J(g)\cap \mathcal{P}}\left\vert [g](x)\cdot \nu
(x)\right\vert {\small d}\mathcal{H}^{N-1}{\small (x)}   
\end{equation}%
for the (minimal) volume swept out by\ separations and switches among
approximations $f_{n}$ that determine $(g,G)$. \ Relation \eqref{formula for volume swept out} provides answers in the setting of Choksi and Fonseca to
the questions raised at the end of Section 1.1: \ $\mathcal{V}^{\left\vert \cdot
\right\vert }(g,G;\mathcal{P})$ has both a bulk disarrangement density $\left\vert
tr(\nabla g-G)\right\vert =\left\vert trM\right\vert $ and an interfacial
disarrangement density $\left\vert [g]\cdot \nu \right\vert $. \ Similarly,
Theorem \ref{th2} shows that the (minimal) volume swept out by separations alone,
${\small \mathcal{V}}^{+}{\small %
(g,G;\mathcal{P})}$, has the bulk disarrangement density $%
(trM)^{+} $ and the interfacial disarrangement density $([g]\cdot \nu )^{+}$, 
with a corresponding result for $\mathcal{V}^{-}(g,G;\mathcal{P})$, the (minimal) volume swept
out by switches and interpenetrations 
(the approximations $f_n$ in the Choksi-Fonseca setting are not required to be injective, so that interpenetrations can arise there, unlike in the setting of Del Piero-Owen).

\subsection{Summary of the research presented in the present article}

In the proof of Theorem \ref{th2} given in \cite{OwenParoni}, the significant part
of the argument addresses the verification of the inequality%
\begin{equation}
H^{\left\vert \cdot \right\vert }(A,B)\leq \left\vert tr(A-B)\right\vert 
\label{upper bound OwenParoni}
\end{equation}%
where $H^{\left\vert \cdot \right\vert }(A,B)$ is given by the right-hand
side of \eqref{relaxed bulk density} with $\psi ([u],\nu _{u})$ replaced by $%
\psi ^{\left\vert \cdot \right\vert }([u],\nu )$ $=\left\vert [u]\cdot \nu
\right\vert $. \ This inequality was proved in \cite{OwenParoni} by
constructing a family $u_{\varepsilon }$ of piecewise affine mappings on the
unit cube $Q$ each of whose jump set $J(u_{\varepsilon })$ is formed by two
(planar) ends and by a lateral surface constructed from solution curves of
the differential equation $\dot{x}=(A-B)x$. \ The lateral surface, by
construction, contributes nothing to the integral $\int_{J(u)\cap Q_{\eta
}}\left\vert [u]\cdot \nu \right\vert \,d\mathcal{H}^{N-1}$, and the
contributions of the two ends can be calculated explicitly for $A-B$ lying
in a dense subset of $\mathbb{R}^{N\times N}$. \ Proposition 5.2 of \cite%
{choksifonseca} provides sufficient regularity of $H^{\left\vert \cdot
\right\vert }(A,B)$ to establish \eqref{upper bound OwenParoni} for all $%
A-B\in \mathbb{R}^{N\times N}$.

As one of the main results in this article, we provide an alternate, shorter
proof of \eqref{upper bound OwenParoni} that employs a different family $%
u_{\varepsilon }$ of piecewise affine mappings that does not involve
solution curves of $\dot{x}=(A-B)x$. Our approach is based on the following
observation. \ With $A,B\in \mathbb{R}^{N\times N}$, $p>1$, and with $%
Q=(-1/2,1/2)^{N}$ there hold 
\begin{equation}\label{dimension N inequalities}
\begin{split}
\left\vert tr(A-B)\right\vert  \leq \inf \bigg\{ & \int_{J(u)}\left\vert [u](x)\cdot \nu (x)\right\vert d%
\mathcal{H}^{N-1}(x)\,:\,u\in SBV(Q;\mathbb{R}^{N}),  \\
&\;{\small u(x)}=Ax\text{ on }\partial Q,\,\nabla u\in L^{p}(Q),\text{
}\,\int_{Q}\nabla u(x)d\mathcal{L}^{N}(x)\,=B\bigg\}   \\
\leq \inf \bigg\{& \int_{J(u)}\left\vert [u](x)\cdot \nu
(x)\right\vert d\mathcal{H}^{N-1}(x)\,:\,u\in SBV(Q;\mathbb{R}^{N}), \\
&\;{\small u(x)}=0\text{ on }\partial Q,\,\nabla u=B-A \; a.e.\bigg\}.
\end{split}
\end{equation}%
\newline
The first follows by moving the absolute value outside the integral and
using the Gauss-Green Theorem for the space $SBV(Q;\mathbb{R}^{N})$ of
special functions of bounded variation, while the second follows by noting
that if $u$ satisfies the last set of conditions, then the function $%
x\longmapsto u(x)+Ax$ satisfies the first set of conditions.{\LARGE \ }\ In
this paper, we wish to show that%
\begin{equation}\label{upper bound inequality}
\begin{split}
\inf \bigg\{ \int_{J(u)}\left\vert [u](x)\cdot \nu (x)\right\vert 
{\small d\mathcal{H}}^{N-1}{\small (x)\,:\,u\in SBV(Q;R}^{N}{\small ),}\;\;\;& \\
{\small u(x)=0}\text{ on }{\small \partial Q,\,\nabla u=B-A}\;a.e.\bigg\} & \leq \left\vert tr(A-B)\right\vert 
\end{split}
\end{equation}%
so that the two infima in \eqref{dimension N inequalities} have common value 
$\left\vert tr(A-B)\right\vert $. \ 

The second main contribution of the present research concerns the
alternative approach to structured deformations and to relaxed energies due
to Ba\'ia, Matias, and Santos \cite{bms} discussed at the end of Subsection \ref{sss}. According to that
discussion the second infimum in \eqref{dimension N inequalities} (see \eqref{relaxed bulk density in BMS} and Remark \ref{bms-remark})
\begin{equation}\label{bms infimum}
\begin{split}
\inf \bigg\{ & \int_{J(u)}\left\vert [u](x)\cdot \nu (x)\right\vert 
 d\mathcal{H}^{N-1}(x)\,:\,u\in SBV(Q;R^{N}), \\
&\; u(x)=0\text{ {\small on} }{\small \partial Q,}
{\small \,\nabla u=B-A}\text{ a.e. in } Q \bigg\}  
\end{split}
\end{equation}%
is the bulk disarrangement density for the same interfacial disarrangement
density $\psi ^{\left\vert \cdot \right\vert }([u],\nu )$ \eqref{psi absolute} studied by Owen and Paroni in the setting of Choksi and Fonseca. \
Consequently, our proof of \eqref{upper bound inequality} establishes the
equality of the bulk disarrangement densities obtained in two different
settings for structured deformations. \ Thus, the geometrical significance
of the expression $\left\vert tr(A-B)\right\vert $ described in \cite%
{OwenParoni}, namely, a volume density of volume swept out by non-smooth,
submacroscopic geometrical changes, is strengthened by the fact that one and
the same expression arises from two different schemes of relaxation. \ We
note that the two different schemes of relaxation also deliver the same formula
for the (relaxed) interfacial disarrangement density $h$: $h=\psi
^{\left\vert \cdot \right\vert }$(see \cite{OwenParoni} for the routine
verification that applies to both schemes).

The explicit formulas for disarrangement densities considered here in the
context of structured deformations will provide scalar fields that can enter
as variables in constitutive relations for the response of three-dimensional
bodies. \ For this purpose, frame-indifferent variants of the specific
fields obtained here are available through known factorizations of
structured deformations in which the factor that tracks disarrangements is
unchanged under changes in frame \cite{del piero owen}. \ Our explicit
formulas also are starting points for the study of examples in other
contexts involving structured deformations: second-order structured
deformations \cite{OwenParoni2} in which second gradients and their limits
enter into submacroscopic changes in geometry, as well as processes for
dimension reduction \cite{MatiasSantos} in the presence of disarrangements
that describe thin structures undergoing submacroscopic slips, separations,
and switches.

In Section 2 we provide a "tilted cube" construction for the family $%
u_{\varepsilon }$ of functions employed in proving \eqref{upper bound inequality}. 
The common orientation of the tilted cubes is determined in
Section 3 by means of a known result on the isotropic vectors of symmetric
linear mappings. The proof of \eqref{upper bound inequality} is completed in
Section 4, and the paper concludes in Section 5 with some additional
explicit formulas for disarrangement densities.

During the review of this article, the research \cite{silhavy} was brought to our attention in which explicit formulas for the bulk and interfacial relaxed energies are established for a broad class of purely interfacial initial energies that includes the ones studied here.

\section{\protect Proof of the upper bound inequality}

In what follows, a proof of \eqref{upper bound inequality} is given. The
proof requires the following instance of Lemma 4.3 in \cite{Matias lemma}.

\begin{lemma}
Let $M\in \mathbb{R}^{N\times N}$ and a bounded open set $\Omega \subset 
\mathbb{R}^{N}$ be given, with $\Omega$ having Lipshitz boundary. \ There exist a number $C(N)>0$, independent of $%
M$ and $\Omega $, and $u\in SBV(\Omega;\mathbb{R}^{N})$ such that

\begin{enumerate}
\item $u|_{\partial \Omega }=0$

\item $\nabla u=M,$ \ \ $\mathcal{L}^{N}-a.e.$ on $\Omega $

\item $\left\vert D^{s}u\right\vert (\Omega )\leq C(N)\left\Vert
M\right\Vert \mathcal{L}^{N}(\Omega ).$
\end{enumerate}
\end{lemma}

Here, $\nabla u$ and $D^{s}u$ denote the absolutely continuous and the
singular parts of the distributional derivative $Du=\nabla u\,\mathcal{L}%
^{N}+D^{s}u$ of $u$, and $\left\vert D^{s}u\right\vert $ denotes the total
variation of the singular part. In addition, \ $\left\Vert M\right\Vert
:=(tr(M^{T}M))^{1/2\text{ }}$ is the Euclidean norm of the matrix $M$. \ We
shall now use the Lemma to verify \eqref{upper bound inequality} for $M=A-B$%
. \ To this end, let an integer $n\geq 1$ be given and consider the frame%
\begin{equation*}
\mathcal{F}_{n}:=Q\setminus \cl{(1-\tfrac{2}{n+2})Q}.
\end{equation*}%
We may apply the Lemma to obtain an $SBV$ function $u^{(n)}:$ $\mathcal{F}_{n}\rightarrow \mathbb{R}^{N}$ satisfying

\begin{itemize}
\item $u^{(n)}|_{\partial \mathcal{F}_{n}}=0$

\item $\nabla u^{(n)}=M$, $\mathcal{L}^{N}-a.e.$ on $\mathcal{F}_{n}$

\item the total variation $\int_{J(u^{(n)})}\left\vert [u^{(n)}]\right\vert
(x)d\mathcal{H}^{N-1}(x)$ of $u^{(n)}$ satisfies%
\begin{equation}
\int_{J(u^{(n)})}\lvert [u^{(n)}]\rvert (x)d\mathcal{H}%
^{N-1}(x)\leq C(N)\left\Vert M\right\Vert \left(1-(1-\tfrac{2}{n+2})^{N}\right)
\label{Fn bound on jumps}
\end{equation}
\end{itemize}

In preparation for defining an appropriate function $u$ on $Q\backslash 
\cl{\mathcal{F}_{n}}=(1-\tfrac{2}{n+2})Q$, we write $\hat{M}:=\tfrac{1}{2}(M+M^{T})$ for
the symmetric part of $M$, and we choose an orthonormal basis $%
e_{i},i=1,\ldots ,N$ of $\mathbb{R}^{N}$ that consists of eigenvectors of $%
\hat{M}$:%
\begin{equation*}
\hat{M}e_{i}=\lambda _{i}e_{i},\ \ i=1,\ldots ,N.
\end{equation*}%
We let $m$ be a positive integer and cover $(1-\tfrac{2}{n+2})Q$ by a
collection $\mathcal{C}_{n,m}$ of congruent, non-overlapping open cubes $%
C_{n,m}^{k}$, $k=1,\ldots ,K_{n,m}$, each of edge-length $1/m$, each with
the $i^{th}$ pair of opposite faces orthogonal to the unit vector $Re_{i}$,
for $i=1,\ldots ,N$. Here, $R$ is an orthogonal $N\times N$ matrix, $%
RR^{T}=R^{T}R=I$, to be determined presently. \ We require in addition that
each each cube $C_{n,m}^{k}$ satisfies 
\begin{equation}
(1-\tfrac{2}{n+2})Q\cap C_{n,m}^{k}\neq \emptyset .
\label{interior condition}
\end{equation}%
\ We denote by $c_{n,m}^{k}$ the center of $C_{n,m}^{k}$, and we define $u_{n,m}:(1-\frac{2}{n+2})Q\to\R{N}$%
\begin{equation}\label{definition of unp}
u_{n,m}(x):=\begin{cases}
M(x-c_{n,m}^{k}) & \text{if $x\in (1-\tfrac{2}{n+2})Q\cap C_{n,m}^{k}$ for some $k=1,\ldots ,K_{n,m}$,} \\ 
0 & \text{if $x\in (1-\tfrac{2}{n+2})Q\backslash \cup_{k=1}^{K_{n,m}} C_{n,m}^{k}$.}%
\end{cases}%
\end{equation}%

Using standard reasoning we conclude that $u_{n,m}\in
SBV((1-\tfrac{2}{n+2})Q;\mathbb{R}^{N})$ with $\nabla u_{n,m}=M$, \thinspace 
$\mathcal{L}^{N}-a.e.$ on $(1-\tfrac{2}{n+2})Q$. \ Moreover, the trace of $%
u_{n,m}$ on $\partial ((1-\tfrac{2}{n+2})Q)$ is bounded pointwise by $\frac{%
\sqrt{N}}{2m}\left\Vert M\right\Vert $. \ Consequently, the function $%
u_{m}^{(n)}:Q\rightarrow \mathbb{R}^{N}$ defined by 
\begin{equation*}
u_{m}^{(n)}(x):=\begin{cases}
u^{(n)}(x) & \text{for $x\in \mathcal{F}_{n}$,} \\ 
u_{n,m}(x) & \text{for $x\in (1-\tfrac{2}{n+2})Q$}
\end{cases}
\end{equation*}
belongs to $SBV(Q;\R{N})$, has gradient $M$, $\mathcal{L}^{N}-$ a.e., and has zero trace on $%
\partial Q$. \ Moreover, the jump set of $u_{m}^{(n)}$ satisfies%
\begin{equation}
J(u_{m}^{(n)})\subset J(u^{(n)})\cup \partial ((1-\tfrac{2}{n+2})Q)\cup
J(u_{n,m}).  \label{pieces of J(u)}
\end{equation}%
Since $u_{m}^{(n)}$ has outer trace $0$ on $\partial (1-\tfrac{2}{n+2})Q,$\
there holds for $\mathcal{H}^{N-1}-a.e.$ $x\ $\ in $\partial ((1-\tfrac{2}{n+2})Q)$%
\begin{equation}
\left\vert \lbrack u_{m}^{(n)}](x)\right\vert \leq \frac{\sqrt{N}}{m}%
\left\Vert M\right\Vert  \label{jump u on frame}
\end{equation}%
and, consequently, 
\begin{equation} \label{estimate on inner square}
\int_{\partial ((1-\tfrac{2}{n+2})Q)}\left\vert [u_{m}^{(n)}](x)\cdot \nu
(x)\right\vert d\mathcal{H}^{N-1}(x) \leq \frac{\sqrt{N}}{m}\left\Vert
M\right\Vert 2N(1-\tfrac{2}{n+2})^{N-1}. 
\end{equation}%
We note from \eqref{Fn bound on jumps} that%
\begin{equation}
\int_{J(u^{(n)})}\left\vert [u_{m}^{(n)}](x)\cdot \nu (x)\right\vert d%
\mathcal{H}^{N-1}(x)\leq C(N)\left\Vert M\right\Vert (1-(1-\tfrac{2}{n+2}%
)^{N}),  \label{estimate on J(un)}
\end{equation}%
and it remains to obtain a corresponding estimate for $\int_{J(u_{n,m})}%
\left\vert [u_{m}^{(n)}](x)\cdot \nu (x)\right\vert d\mathcal{H}^{N-1}(x)$.
\ To this end, we note that%
\begin{equation}
J(u_{n,m})\subset \bigcup _{k=1}^{K_{n,m}}\partial C_{n,m}^{k}\,,
\label{jump set of unp}
\end{equation}%
and we shall seek an upper bound for $\int_{\cup _{k=1}^{K_{n,m}}\partial
C_{n,m}^{k}}\left\vert [u_{m}^{(n)}](x)\cdot \nu (x)\right\vert d\mathcal{H}%
^{N-1}(x)$. \ For each $k=1,\ldots ,K_{n,m}$ and $i=1,\ldots ,N\allowbreak ,$
we denote by $\phi _{n,m}^{k,i+}$ and $\phi _{n,m}^{k,i-}$ the two faces of
the cube $C_{n,m}^{k}\in $ $\mathcal{C}_{n,m}$ orthogonal to $R\,e_{i}$. We
note that one face $\phi _{n,m}^{k,i+}$of $C_{n,m}^{k}$ has outer normal $%
\nu _{i,k}^{+}=$ $+R\,e_{i}$, while the opposite face $\phi
_{n,m}^{k,i-}$ has outer normal $\nu _{i,k^{{}}}^{-}=$ $-R\,e_{i}$.

\ We suppose now that the face $\phi _{n,m}^{k,i+}$of $C_{n,m}^{k}\in $ $%
\mathcal{C}_{n,m}$ satisfies%
\begin{equation}
\phi _{n,m}^{k,i+}\subset (1-\tfrac{2}{n+2})Q.  \label{inside sides plus}
\end{equation}%
Then there is a cube $C_{n,m}^{k^{\prime }}\in \mathcal{C}_{n,m}$ that
shares the given face with $C_{n,m}^{k}$, and we have at each point $x\in $ $%
\phi _{n,m}^{k,i+}$%
\begin{eqnarray*}
\lbrack u_{m}^{(n)}](x)\cdot \nu (x) &=&(M(x-c_{n,m}^{k_{{}}^{\prime
}})-M(x-c_{n,m}^{k_{{}}}))\cdot \nu _{i,k}^{+}(x) \\
&=&M(c_{n,m}^{k_{{}}}-c_{n,m}^{k_{{}}^{\prime }})\cdot \nu _{i,k}^{+}(x) \\
&=&M(-\tfrac{1}{m}\nu _{i,k}^{+}(x))\cdot \nu _{i,k}^{+}(x) \\
&=&-\frac{1}{m}\hat{M}R\,e_{i}\cdot R\,e_{i}
\end{eqnarray*}%
so that 
\begin{equation}\label{contribution of plus side}
\begin{split}
\int_{\phi _{n,m}^{k,i+}}\left\vert [u_{m}^{(n)}](x)\cdot \nu (x)\right\vert
d\mathcal{H}^{N-1}(x) =&\int_{\phi _{n,m}^{k,i+}}\frac{1}{m}\vert 
\hat{M}R\,e_{i}\cdot R\,e_{i}\vert d\mathcal{H}^{N-1}(x)   \\
=&\frac{1}{m^{N}}\vert \hat{M}R\,e_{i}\cdot R\,e_{i}\vert .
\end{split}
\end{equation}%
The same argument shows that if 
\begin{equation}
\phi _{n,m}^{k,i-}\subset (1-\tfrac{2}{n+2})Q  \label{inside sides minus}
\end{equation}%
then%
\begin{equation}
\int_{\phi _{n,m}^{k,i-}}\left\vert [u_{m}^{(n)}](x)\cdot \nu (x)\right\vert
d\mathcal{H}^{N-1}(x)=\frac{1}{m^{N}}\vert \hat{M}R\,e_{i}\cdot
R\,e_{i}\vert .  \label{contribution of minus}
\end{equation}%
If \eqref{inside sides plus} holds for $i=1,\ldots ,N$, then we may sum the
last relation over $i$ to conclude that%
\begin{equation} \label{basic trace inequality}
\begin{split}
\sum_{i=1}^{N}\int_{\phi _{n,m}^{k,i+}}\left\vert [u_{m}^{(n)}](x)\cdot \nu
(x)\right\vert d\mathcal{H}^{N-1}(x) =&\frac{1}{m^{N}}\sum_{i=1}^{N}
\vert \hat{M}R\,e_{i}\cdot R\,e_{i}\vert   \\
\geq &\frac{1}{m^{N}}\left\vert \sum_{i=1}^{N}\hat{M}R\,e_{i}\cdot
R\,e_{i}\right\vert   \\
=&\frac{1}{m^{N}}\left\vert \sum_{i=1}^{N}R^{T}\hat{M}R\,e_{i}\cdot
\,e_{i}\right\vert   \\
=&\frac{1}{m^{N}}\vert tr(R^{T}\hat{M}R)\vert =\frac{1}{m^{N}}%
\left\vert trM\right\vert.
\end{split}
\end{equation}%
In \eqref{basic trace inequality} equality holds if and only if all of the
numbers $\hat{M}R\,e_{i}\cdot R\,e_{i}\,$, $i=1,\ldots ,N$, have the same
sign: 
\begin{equation}
(\hat{M}R\,e_{i}\cdot R\,e_{i}\,)(\hat{M}R\,e_{j}\cdot R\,e_{j}\,)\geq 0,%
\text{ for }i,j=1,\ldots ,N.  \label{sign inequality}
\end{equation}%
The last two inequalities lead us to consider the problem%
\begin{equation}
\text{Find }\min_{RR^{T}=I}\sum_{i=1}^{N}\vert \hat{M}R\,e_{i}\cdot
R\,e_{i}\vert \geq \vert tr\,\hat{M}\vert =\left\vert
trM\right\vert ,  \label{minimum condition for R}
\end{equation}%
with equality holding if and only if there exists an orthogonal matrix $R$
satisfying \eqref{sign inequality}.

\section{Aside on isotropic vectors}

We note that the sign inequality \eqref{sign inequality} suggests looking
for unit vectors $v$ such that 
\begin{equation}
\hat{M}v\cdot v=0,  \label{definition isotropic}
\end{equation}%
the \textit{isotropic vectors} for $\hat{M}$ \cite{gatech}. \ In particular,
in the special case $tr\,\hat{M}=0$, the existence of $N$ mutually
orthogonal isotropic vectors $v_{1},\ldots ,v_{N}$ would insure that the
matrix $R$ defined by $R\,e_{i}=v_{i}$ for $i=1,\ldots ,N$ would satisfy \eqref{minimum condition for R} in the form $0=0$. \ More generally, even when \ $%
tr\,\hat{M}\neq 0$, the existence of isotropic vectors is useful. \ In
fact, the symmetric matrix $\hat{M}-\frac{1}{N}(tr\,\hat{M})I$ has zero
trace, so we suppose that there exist $N$ mutually orthogonal \textit{%
isotropic} \textit{unit} vectors $v_{1},\ldots ,v_{N}$ for $\hat{M}-\frac{1}{%
N}(tr\,\hat{M})I$. \ The relation \eqref{definition isotropic} with $\hat{M}$
replaced by $\hat{M}-\frac{1}{N}(tr\,\hat{M})I$ then becomes%
\begin{equation*}
0 =\Big(\hat{M}-\frac{1}{N}(tr\,\hat{M})I\Big)v_{i}\cdot v_{i} =\hat{M}v_{i}\cdot v_{i}-\frac{tr\,\hat{M}}{N}
\end{equation*}%
so that $\hat{M}v_{i}\cdot v_{i}=\frac{tr\,\hat{M}}{N}$ for $i=1,\ldots ,N$.
Again, if we define a linear mapping $R$ on $\mathbb{R}^{N}$by $%
R\,e_{i}=v_{i}$ for $i=1,\ldots ,N$ then $R$ is orthogonal, it satisfies the
sign inequality for $\hat{M}$ \eqref{sign inequality}, and it delivers
equality in \eqref{minimum condition for R} in the form $\sum_{i=1}^{N}\left%
\vert \frac{tr\,\hat{M}}{N}\right\vert =\vert tr\,\hat{M}\vert $.

The following result (\cite{gatech}, Corollary 15) provides the desired
existence of complete orthonormal sets of isotropic vectors.{}

\begin{theorem}
A symmetric matrix $A\in R^{N\times N}$ possesses an orthonormal set of $N$
isotropic vectors if and only if $trA=0$.
\end{theorem}

This theorem and the preceding discussion permit us to conclude: for every
matrix $M\in \mathbb{R}^{N\times N},$ 
\begin{equation}\label{solution to minimum}
\begin{split}
\min_{RR^{T}=I}\sum_{i=1}^{N}\vert MR\,e_{i}\cdot R\,e_{i}\vert
=&\min_{RR^{T}=I}\sum_{i=1}^{N}\vert \hat{M}R\,e_{i}\cdot
R\,e_{i}\vert   \\
=&\vert tr\,\hat{M}\vert =\left\vert trM\right\vert,
\end{split}
\end{equation}%
and a minimizing rotation matrix $R$ is one carrying the orthonormal basis
of $\mathbb{R}^{N}$ consisting of eigenvectors of $\hat{M}$ into an
orthonormal basis of $\mathbb{R}^{N}$ consisting of isotropic vectors of $\hat{M%
}-\frac{1}{N}(tr\,\hat{M})I$. For this minimizing rotation matrix, we have%
\begin{equation}
\vert \hat MR\,e_{i}\cdot R\,e_{i}\vert =\frac{1}{N}\left\vert
trM\right\vert \text{ for }i=1,\ldots ,N\text{.}
\label{ith term in minimizing sum}
\end{equation}

We remark that minimizers are not unique, in general, even when one
eliminates trivial permutations of isotropic vectors. \ In fact, for $N=3$
there are examples of minimizers for which two of the three terms in $%
\sum_{i=1}^{3}\vert \hat{M}R\,e_{i}\cdot R\,e_{i}\vert $ vanish,
while the third equals $\left\vert trM\right\vert $, so that only two of the
three vectors $R\,e_{i}$ are isotropic vectors for $\hat{M}$ .

For the convenience of the reader, we provide the recursive step used in
proving the existence of orthonormal bases made up of isotropic vectors for
a traceless symmetric matrix $A\in \mathbb{R}^{N\times N}$. \ We interpret $%
A $ in the usual way as a linear mapping on $\mathbb{R}^{N},$ endowed with
the standard inner product. \ Then the nullspace $KerA$ of $A$ and its
orthogonal complement $(KerA)^{\perp }$ are complementary $A$-invariant
subspaces of $\mathbb{R}^{N}$, and all vectors in $KerA$ are isotropic
vectors for $A$. \ If \ $(KerA)^{\perp }$ is the zero subspace, then $A=0$
and every vector in $\mathbb{R}^{N}$ is an isotropic vector for $A$, and
every orthonormal basis of $\mathbb{R}^{N}$ meets the desired requirement. \
If $(KerA)^{\perp }$ is not the zero subspace, then we seek additional
isotropic vectors for $A$ in $(KerA)^{\perp }$. \ To this end, the traceless
symmetric linear mapping $A\neq 0$ has both positive and negative
eigenvalues so that 
\begin{equation*}
\min_{\left\vert u\right\vert =1}Au\cdot u<0<\max_{\left\vert u\right\vert
=1}Au\cdot u
\end{equation*}%
and, since the unit sphere in $\mathbb{R}^{N}$ is connected and since the
quadratic form $u\longmapsto Au\cdot u$ is continuous, there exists a unit
vector $v_{1}\in \mathbb{R}^{N}$ such that $Av_{1}\cdot v_{1}=0$. \ Writing $%
v_{1}$ as a sum of two orthogonal vectors, one in $KerA$ and the other in $%
(KerA)^{\perp }$ and using the invariance of $(KerA)^{\perp }$ under $A$
shows that we may without loss of generality assume that $v_{1}\in
(KerA)^{\perp }$. \ The linear span $Lsp(KerA\cup \left\{ v_{1}\right\} )$
has dimension one larger than that of $KerA$ and consists solely of
isotropic vectors for $A$. \ Consequently, we need to search for isotropic
vectors of $A$ in $(Lsp(KerA\cup \left\{ v_{1}\right\} ))^{\perp }$ which
has dimension one less than $(KerA)^{\perp }$. \ To procede further, we
define a linear mapping $A_{1}$ on $\mathbb{R}^{N}$ by%
\begin{equation}
A_{1}=A-v_{1}\otimes Av_{1}-Av_{1}\otimes v_{1}  \label{definition of A1}
\end{equation}%
where the formula $(a\otimes b)v:=(b\cdot v)a$, for all $a,b,v\in \mathbb{R}%
^{N}$, defines the standard tensor product $a\otimes b\in Lin(\mathbb{R}^{N};\mathbb{R}^{N})$. \ From the fact that $v_{1}$ is an isotropic vector for $A$
and from the formula $tr((a\otimes b)=a\cdot b$ it is easy to see that $%
A_{1} $ is traceless; because $(a\otimes b)^{T}=b\otimes a$, it follows that 
$A_{1} $ is symmetric. In addition, if $v\in (Lsp(KerA\cup \left\{
v_{1}\right\} ))^{\perp }$ is an isotropic vector for $A_{1}$, then we have
not only $v\cdot v_{1}=0$ but also 
\begin{eqnarray*}
0 &=&A_{1}v\cdot v \\
&=&(Av-(Av_{1}\cdot v)v_{1}-(v_{1}\cdot v)Av_{1})\cdot v \\
&=&Av\cdot v-(Av_{1}\cdot v)(v_{1}\cdot v)-(v_{1}\cdot v)(Av_{1}\cdot v) \\
&=&Av\cdot v\text{.}
\end{eqnarray*}%
Thus, every isotropic vector for $A_{1}$ that is in $(Lsp(KerA\cup \left\{
v_{1}\right\} ))^{\perp }$ is an isotropic vector for $A$, and $\dim
((Lsp(KerA\cup \left\{ v_{1}\right\} ))^{\perp })=\dim ((KerA)^{\perp })-1$.
\ To be able to apply the forgoing considerations to $A_{1}$, we need only
show that $(Lsp(KerA\cup \left\{ v_{1}\right\} ))^{\perp }$ \ is invariant
under $A_{1}$. \ To this end, let $v\in (Lsp(KerA\cup \left\{ v_{1}\right\}
))^{\perp }$, $v_{\kappa }\in $\ $KerA$, and $\alpha \in \mathbb{R}$ be
given, and consider%
\begin{eqnarray*}
A_{1}v\cdot (v_{\kappa }+\alpha v_{1}) &=&A_{1}v\cdot v_{\kappa
}+A_{1}v\cdot \alpha v_{1} \\
&=&v\cdot A_{1}v_{\kappa }+\alpha v\cdot A_{1}v_{1} \\
&=&0+\alpha v\cdot (Av_{1}-(Av_{1}\otimes v_{1})v_{1}-(v_{1}\otimes
Av_{1})v_{1}) \\
&=&\alpha v\cdot (Av_{1}-(v_{1}\cdot v_{1})Av_{1}-(Av_{1}\cdot v_{1})v_{1})
\\
&=&\alpha v\cdot (Av_{1}-Av_{1}-0)=0.
\end{eqnarray*}%
We may conclude that $A_{1}v\in (Lsp(KerA\cup \left\{ v_{1}\right\}
))^{\perp }$ as desired. In the third line of the above
computation we have used the side-calculation
\begin{eqnarray*}
v\cdot A_{1}v_{\kappa } &=
&v\cdot (A-v_{1}\otimes Av_{1}-Av_{1}\otimes v_{1})v_{\kappa } \\
&=&v\cdot Av_{\kappa }-(Av_{1}\cdot v_{\kappa })(v\cdot v_{1})-(v_{1}\cdot
v_{\kappa })(v\cdot Av_{1})=0.
\end{eqnarray*}%
The first term on the last line vanishes because $v_{\kappa }\in $%
{\large \ }$KerA$, the second vanishes because $v\in
(Lsp(KerA\cup \left\{ v_{1}\right\} ))^{\perp }$ and the third
vanishes because $v_{1}\in (KerA)^{\perp }${\large .} \ The search for
isotropic vectors for $A$ on the $A$-invariant subspace $(KerA)^{\perp }$
may now be replaced by the search for isotropic vectors for $A_{1}$ on the $%
A_{1}$-invariant subspace $(Lsp(KerA\cup \left\{ v_{1}\right\} ))^{\perp }$
of dimension one less than that of $(KerA)^{\perp }$.\ 

\section{Completion of the proof of the upper bound inequality}

We may use \eqref{ith term in minimizing sum} and the formulas \eqref{contribution of plus side}, \eqref{contribution of minus} to conclude
that: if $C_{n,m}^{k}$ has a face $\phi _{n,m}^{k,i\pm }\subset (1-\frac{2}{%
N+2})Q$, then 
\begin{equation}
\int_{\phi _{n,m}^{k,i\pm }}\left\vert [u_{m}^{(n)}](x)\cdot \nu
(x)\right\vert d\mathcal{H}^{N-1}(x)=\frac{\left\vert trM\right\vert }{Nm^{N}%
}=\frac{\left\vert trM\right\vert }{N}\mathcal{L}^{N}(C_{n,m}^{k}).
\label{basic trace relation}
\end{equation}%
On the other hand, if a face $\phi _{n,m}^{k,i\pm }$ of $C_{n,m}^{k}\in 
\mathcal{C}_{n,m}$ fails to satisfy $\phi _{n,m}^{k,i\pm }\subset (1-\frac{2%
}{N+2})Q$, then the argument used to verify \eqref{basic trace relation} may
be applied to $\phi _{n,m}^{k,i\pm }\cap $ $(1-\tfrac{2}{n+2})Q$ $\ $to
conclude that%
\begin{equation}
\int_{\phi _{n,m}^{k,i\pm }\cap (1-\tfrac{2}{n+2})Q}\left\vert
[u_{m}^{(n)}](x)\cdot \nu (x)\right\vert d\mathcal{H}^{N-1}(x)\leq \frac{%
\left\vert trM\right\vert }{N}\mathcal{L}^{N}(C_{n,m}^{k}).
\label{bound for contributions of sides}
\end{equation}

We now consider the cube $C_{n,m}^{1}\in \mathcal{C}_{n,m}$ and choose $%
V_{n,m}^{1}$, one of its $2^{N}$ vertices. \ Exactly $N$ faces $\phi
_{{}}^{1,j},j=1,\ldots ,N$, of $C_{n,m}^{{1}}$ meet at $V_{n,m}^{1}$. \
Because each cube $C_{n,m}^{k}\in \mathcal{C}_{n,m}$ for $k=1,\ldots ,K_{n,m}$ can be obtained from $C_{n,m}^{1}$ by a unique translation $T_{k}$, the choices $C_{n,m}^{1}$ and $V_{n,m}^{1}$ induce via $T_{k}$ an
assignment of $N$ faces $\phi^{k,j},$ $j=1,\ldots ,N$ $\ $to $%
C_{n,m}^{k}$. \ It is easy to show that for all $k,k^{\prime }=1,\ldots,K_{n,m}$ 
\begin{equation*}
k^{\prime }\neq k\Longrightarrow \big\{ \phi^{k^{\prime},j}:j=1,\ldots,N\big\} \cap \big\{ \phi^{k,j}:j=1,\ldots,N\big\} =\emptyset ,
\end{equation*}%
i.e., the set of $N$ faces assigned to different cubes are disjoint. If we
now apply the mapping%
\begin{equation*}
C_{n,m}^{k}\longmapsto \left\{ \phi^{k,j}:j=1,\ldots ,N\right\} \,
\end{equation*}%
to each cube in the collection 
\begin{equation*}
\mathcal{C}_{n,m}^{\itr}:=\bigg\{ C_{n,m}^{k}\in \mathcal{C}%
_{n,m}:C_{n,m}^{k}\subset (1-\frac{2}{N+2})Q\bigg\},
\end{equation*}%
then all of the faces $\phi^{k,j}$ so obtained will be included
in $(1-\frac{2}{N+2})Q$, and we may apply \eqref{basic trace relation} to
each such face to obtain for each $C_{n,m}^{k}\in \mathcal{C}_{n,m}^{\itr}$ 
\begin{equation*}
\sum_{j=1}^{N}\int_{\phi^{k,j}}\left\vert [u_{m}^{(n)}](x)\cdot \nu
(x)\right\vert d\mathcal{H}^{N-1}(x) =N\frac{\left\vert trM\right\vert }{N}
\mathcal{L}^{N}(C_{n,m}^{k})=\left\vert trM\right\vert \mathcal{L}^{N}(C_{n,m}^{k}).
\label{sum relation on interior cube}
\end{equation*}%
We may sum both sides over the cubes $C_{n,m}^{k}\in \mathcal{C}_{n,m}^{\itr}$
to obtain%
\begin{equation}
\sum_{C_{n,m}^{k}\in \mathcal{C}_{n,m}^{\itr}}\sum_{j=1}^{N}\int_{\phi^{k,j}}\left\vert [u_{m}^{(n)}](x)\cdot \nu (x)\right\vert d\mathcal{H}%
^{N-1}(x)=\left\vert trM\right\vert \mathcal{L}^{N}(\cup _{C_{n,m}^{k}\in 
\mathcal{C}_{n,m}^{\itr}}C_{n,m}^{k})\text{.}
\label{sum relation for all interior cubes}
\end{equation}%
The faces represented on the left hand side need not include all of $%
J(u_{n,m})\subset \cup _{k=1}^{K_{n,m}}\partial C_{n,m}^{k}\,$,\ because
some faces of cubes $C_{n,m}^{k}\in \mathcal{C}_{n,m}^{\itr}$ that are also
faces of cubes $C_{n,m}^{k^{\prime }}\in \mathcal{C}_{n,m}\backslash 
\mathcal{C}_{n,m}^{\itr}$ are left out, while proper subsets $\phi
_{n,m}^{k,i\pm }\cap $ $(1-\tfrac{2}{n+2})Q$ of faces $\phi _{n,m}^{k,i\pm }$
also are left out. \ However, for those parts of $J(u_{n,m})$, we may use \eqref{basic trace relation} and \eqref{bound for contributions of sides} to
estimate the integrals $\displaystyle\int_{\phi _{n,m}^{k,i\pm }\cap (1-\tfrac{2}{n+2}%
)Q}\left\vert [u_{m}^{(n)}](x)\cdot \nu (x)\right\vert d\mathcal{H}^{N-1}(x)$%
, along with the fact that the cubes whose faces contain these parts of $%
J(u_{n,m})$ all must contain points of $\partial (1-\tfrac{2}{n+2})Q$ and
must together cover $\partial (1-\tfrac{2}{n+2})Q$. \ Combining all of these
contributions to $\displaystyle\int_{J(u_{n,m})}\left\vert [u_{m}^{(n)}](x)\cdot \nu
(x)\right\vert d\mathcal{H}^{N-1}(x)$ we obtain%
\begin{equation}\label{final estimate for difference}
\begin{split}
0 \leq &\int_{J(u_{n,m})}\left\vert [u_{m}^{(n)}](x)\cdot \nu
(x)\right\vert d\mathcal{H}^{N-1}(x)-\left\vert trM\right\vert \mathcal{L}%
^{N}(\cup _{C_{n,m}^{k}\in \mathcal{C}_{n,m}^{\itr}}C_{n,m}^{k})   \\
\leq &2\left\vert trM\right\vert \mathcal{L}^{N}(\cup _{C_{n,m}^{k}\in (%
\mathcal{C}_{n,m}\backslash \mathcal{C}_{n,m}^{\itr})}C_{n,m}^{k}).
\end{split}
\end{equation}%
The factor of $2=\frac{2N}{N}$ in the last expression reflects the fact that
the $\mathcal{L}^{N}$-measure of some of the cubes in the collection $\mathcal{C%
}_{n,m}\backslash \mathcal{C}_{n,m}^{\itr}$ has been counted more than once
but no more than $2N$ times through the use of the bound \eqref{bound for contributions of sides}. 
The relations \eqref{final estimate for difference}, \eqref{pieces of J(u)}, \eqref{estimate on inner square}, and \eqref{estimate on J(un)} now yield the relation%
\begin{equation}\label{J(u) master estimate}
\begin{split}
0 \leq &\int_{J(u_{m}^{(n)})}\left\vert [u_{m}^{(n)}](x)\cdot \nu
(x)\right\vert d\mathcal{H}^{N-1}(x)-\left\vert trM\right\vert \mathcal{L}%
^{N}(\cup _{C_{n,m}^{k}\in \mathcal{C}_{n,m}^{\itr}}C_{n,m}^{k})   \\
\leq &{\small 2}\left\vert tr{\small M}\right\vert \mathcal{L}^{N}{\small %
(\cup _{C_{n,m}^{k}\in (\mathcal{C}_{n,m}\backslash \mathcal{C}%
_{n,m}^{\itr})}C_{n,m}^{k})}+\frac{\sqrt{N}}{m}\left\Vert {\small M}%
\right\Vert {\small 2N(1-}\tfrac{2}{n+2}{\small )}^{N-1}   \\
&+{\small C(N)}\left\Vert {\small M}\right\Vert {\small (1-(1-}\tfrac{2}{n+2%
}{\small )}^{N}{\small ).}  
\end{split}
\end{equation}%
We use in turn \eqref{J(u) master estimate} to obtain an upper bound for%
\begin{equation*}
\int_{J(u_{m}^{(n)})}\left\vert [u_{m}^{(n)}](x)\cdot \nu (x)\right\vert d%
\mathcal{H}^{N-1}(x).
\end{equation*}%
Let $\varepsilon >0$ $\ $be given and choose $n$ so large that ${\small C(N)}%
\left\Vert {\small M}\right\Vert {\small (1-(1-}\tfrac{2}{n+2}{\small )}^{N}%
{\small )}<\varepsilon $ and, for such an $n$, choose $m$ so large that \ $%
\tfrac{\sqrt{N}}{m}\left\Vert {\small M}\right\Vert {\small 2N(1-}\tfrac{2}{%
n+2}{\small )}^{N-1}<\varepsilon $. \ Because $(1-\tfrac{2}{n+2})Q$ has
finite $\mathcal{L}^{N}-$measure, we may choose $m$ larger if necessary so
that the cover $\mathcal{C}_{n,m}$ of $(1-\tfrac{2}{n+2})Q$ satisfies $%
\mathcal{L}^{N}(\cup _{C_{n,m}^{k}\in \mathcal{C}_{n,m}}C_{n,m}^{k})<%
\mathcal{L}^{N}((1-\tfrac{2}{n+2})Q)+\varepsilon <1+\varepsilon $. \
Finally, because $\partial (1-\tfrac{2}{n+2})Q$ has zero $\mathcal{L}^{N}-$%
measure and is covered by $\mathcal{C}_{n,m}{\small \backslash }\mathcal{C}%
_{n,m}^{\itr}$, we may again choose $m$ larger, if necessary, so that $%
{\small 2}\left\vert tr{\small M}\right\vert \mathcal{L}^{N}{\small (\cup
_{C_{n,m}^{k}\in (\mathcal{C}_{n,m}\backslash \mathcal{C}%
_{n,m}^{\itr})}C_{n,m}^{k})}<\varepsilon $. \ \ We conclude that for $n$ and $%
m$ so chosen 
\begin{equation*}
\int_{J(u_{m}^{(n)})}\left\vert [u_{m}^{(n)}](x)\cdot \nu (x)\right\vert d%
\mathcal{H}^{N-1}(x) <\left\vert trM\right\vert (1+\varepsilon
)+3\varepsilon=\left\vert trM\right\vert +(\left\vert trM\right\vert +3)\varepsilon
\end{equation*}%
and, since $\varepsilon >0$ was arbitrary, that \eqref{upper bound inequality} holds. \ \ $\square $

\section{ Additional explicit formulas for disarrangement densities}

Our discussion above shows that the particular choice of interfacial measure
of disarrangements 
\begin{equation}
\int_{J(u)\cap \Omega }\left\vert [u]\cdot \nu \right\vert d\mathcal{H}^{N-1}
\label{initial interfacial energy u dot nu}
\end{equation}%
for deformations $u$ of a region $\Omega \subset \mathbb{R}^{N}$ leads in
both the Choksi-Fonseca relaxation scheme \cite{choksifonseca} and in the
Ba\'ia-Matias-Santos relaxation scheme \cite{bms} to one and the same bulk
density of disarrangements%
\begin{equation}
\int_{\Omega }\left\vert tr(\nabla g-G)\right\vert d\mathcal{L}^{N}
\label{relaxed bulk energy u dot nu}
\end{equation}%
for structured deformations $(g,G)$ of that region. \ Moreover, our analysis
here provides an alternative to the proof of this result given in \cite%
{OwenParoni} . \ In that article, it was observed that replacement of $%
\left\vert [u]\cdot \nu \right\vert $ by its positive part $([u]\cdot \nu
)^{+}=\frac{1}{2}(\left\vert [u]\cdot \nu \right\vert +$ $[u]\cdot \nu )$
results in the replacement of $\left\vert tr(\nabla g-G)\right\vert $ by
its positive part $(tr(\nabla g-G))^{+}=\frac{1}{2}(\left\vert tr(\nabla
g-G)\right\vert +$ $tr(\nabla g-G))$ in the relaxed bulk disarrangement
density. \ (An analogous result holds for the negative parts, obtained by
replacing "$+$" by "$-$" in the definition of the positive parts.) \ As
pointed out in \cite{OwenParoni}, $(tr(\nabla g-G)(x))^{+}$ may now be
interpreted as the minimum volume fraction at a point $x\in \Omega $ that
can be swept out by submacroscopic separations associated with deformations $%
u_{n}$ approximating the structured deformation $(g,G)$. \ Moreover, $%
(tr(\nabla g-G)(x))^{-}$ is the minimum volume fraction at $x$ swept out by
submacroscopic switches and interpenetrations, so that $\left\vert tr(\nabla
g-G)(x)\right\vert =(tr(\nabla g-G)(x))^{+}+(tr(\nabla g-G)(x))^{-}$ is the
minimum volume fraction swept out by submacroscopic separations, switches,
and interpenetrations.

The presence of the inner-product $[u]\cdot \nu $ in the initial interfacial
density \eqref{initial interfacial energy u dot nu} tells us that only
normal components of jumps will contribute and that alternative initial
interfacial densities are required in order to capture contributions of
tangential components of jumps. \ In the remainder of this section we shall
provide alternative initial interfacial densities that not only capture
contributions of tangential components of jumps but\ also lead to specific
formulas for the relaxed bulk disarrangement density via the "tilted cube"
construction provided in Sections 2 and 4 above.

Let $a\in \mathbb{R}^{N}$ be given and consider the following replacement
for \eqref{initial interfacial energy u dot nu}
\begin{equation}
\int_{J(u)\cap \Omega }\left\vert [u]\cdot a\right\vert d\mathcal{H}^{N-1}
\label{initial interfacial energy u dot a}
\end{equation}%
in which the normal component $[u]\cdot \nu $ of the jump in $u$ is replaced
by the component $[u]\cdot a$ in the direction of $a$. To follow again the
relaxation scheme in \cite{bms} we let $A,B\in \mathbb{R}^{N\times N}$ be
given and require not only $u\in SBV(Q,\mathbb{R}^{N})$ but also%
\begin{equation}
u|_{\partial Q}=0\text{ \ \ , \ }\nabla u=B-A\text{, \ }\;\mathcal{L}%
^{N}-a.e.\text{ in }Q.  \label{cube conditions}
\end{equation}%
We now may use the Gauss-Green formula and \eqref{cube conditions} to write%
\begin{equation}\label{lower bound for u dot a}
\begin{split}
\int_{J(u)\cap Q}\left\vert [u]\cdot a\right\vert d\mathcal{H}^{N-1}
=&\int_{J(u)\cap Q}\left\vert ([u\cdot a])\nu \right\vert d\mathcal{H}^{N-1} \\
\geq &\left\vert \int_{J(u)\cap Q}([u\cdot a])\nu d\mathcal{H}%
^{N-1}\right\vert   \\
=&\left\vert -\int_{Q}\nabla (u\cdot a)d\mathcal{L}^{N}+\int_{\partial
Q}(u\cdot a)\nu d\mathcal{H}^{N-1}\right\vert   \\
=&\left\vert -\int_{Q}(\nabla u)^{T}a\,d\mathcal{L}^{N}+\int_{\partial
Q}(0\cdot a)\nu d\mathcal{H}^{N-1}\right\vert   \\
=&\left\vert (B-A)^{T}a\right\vert  
\end{split}
\end{equation}%
For the "tilted-cube" construction provided in Sections 2 and 4, we replace
the matrix $M$ by $B-A$ , and the relation \eqref{contribution of plus side}
has here the following counterpart%
\begin{equation*}
\begin{split}
\int_{\phi _{n,m}^{k,i+}}\left\vert [u_{m}^{(n)}](x)\cdot a\right\vert d%
\mathcal{H}^{N-1}(x) =&\int_{\phi _{n,m}^{k,i+}}\left\vert
([u_{m}^{(n)}](x)\cdot a)\nu (x)\right\vert d\mathcal{H}^{N-1}(x) \\
=&\int_{\phi _{n,m}^{k,i+}}\frac{1}{m}\left\vert ((B-A)R\,e_{i}\cdot
a)R\,e_{i}\right\vert d\mathcal{H}^{N-1}(x) \\
=& \frac{1}{m^{N}}\left\vert (R\,e_{i}\cdot
(B-A)^{T}a)\,R\,e_{i}\right\vert ,
\end{split}
\end{equation*}%
and this formula leads to the following counterpart of \eqref{basic trace inequality}:
\begin{equation}\label{basic trace inequality u dot a}
\begin{split}
\sum_{i=1}^{N}\int_{\phi _{n,m}^{k,i+}}\left\vert [u_{m}^{(n)}](x)\cdot
a\right\vert d\mathcal{H}^{N-1}(x) =&\frac{1}{m^{N}}\sum_{i=1}^{N}\left%
\vert (R\,e_{i}\cdot (B-A)^{T}a)\,R\,e_{i}\right\vert   \\
\geq &\frac{1}{m^{N}}\left\vert \sum_{i=1}^{N}(R\,e_{i}\cdot
(B-A)^{T}a)\,R\,e_{i}\right\vert   \\
=&\frac{1}{m^{N}}\left\vert (B-A)^{T}a\right\vert.
\end{split}
\end{equation}%
The method employed in Sections 2 and 4 (where the symbol $M$ was used in
place of $B-A$) then requires the choice of a rotation $R$ for which
equality holds in the second line of \eqref{basic trace inequality u dot a}.
\ If $(B-A)^{T}a\neq 0$ we may choose $R$ to be any rotation satisfying $%
R\,e_{1}=$ $(B-A)^{T}a\,/\,\left\vert (B-A)^{T}a\right\vert $, and this
requirement is then met, because $(R\,e_{i}\cdot (B-A)^{T}a)\,R\,e_{i}=0$
for $i=2,\ldots ,N$. \ If $(B-A)^{T}a=0$, then $R$ can be chosen
arbitrarily, for example, $R=I$ suffices.

These observations show that the analysis in Section 4 for \eqref{initial interfacial energy u dot nu} may be carried out step by step for the
alternative initial density \eqref{initial interfacial energy u dot a},
provided that we replace everywhere in Section 4 $\left\vert trM\right\vert
=\left\vert tr(B-A)\right\vert $ by $\,\left\vert (B-A)^{T}a\right\vert $,
the Euclidean norm of the vector $(B-A)^{T}a$. If we now define%
\begin{equation}\label{definition of H(A,B,a)}
\begin{split}
H(A,B,a) :=\inf \bigg\{ \int_{J(u)} & \left\vert [u](x)\cdot
a\right\vert d\mathcal{H}^{N-1}(x)\,:\,u\in SBV(Q;\mathbb{R}^{N}),   \\
& \; u \mid_{\partial Q}=0\text{, }\nabla u=B-A\text{ }a.e.\bigg\},
\end{split}
\end{equation}%
then our observations amount to the formula%
\begin{equation}
H(A,B,a)=\left\vert (B-A)^{T}a\right\vert  \label{formula for H(A,B,a)}
\end{equation}%
for the relaxed bulk energy density corresponding to the initial interfacial
energy \eqref{initial interfacial energy u dot a} and arising from the
scheme \cite{bms}. Moreover, an argument similar to that used in
establishing \eqref{dimension N inequalities} shows that the formula \eqref{formula for H(A,B,a)} also holds for the relaxed bulk disarrangement
density according to \cite{choksifonseca}. In the context of a given
structured deformation $(g,G)$ on a region $\Omega $, \eqref{formula for H(A,B,a)} implies that the particular choice of initial interfacial
disarrangement 
\begin{equation}
\int_{J(u)\cap \Omega }\left\vert [u]\cdot a\right\vert d\mathcal{H}^{N-1}
\label{initial bulk energy u dot a bis}
\end{equation}%
for deformations $u$ of a region $\Omega \subset \mathbb{R}^{N}$ leads in
both the Choksi-Fonseca relaxation scheme \cite{choksifonseca} and in the
Ba\'ia-Matias-Santos relaxation scheme \cite{bms} to one and the same relaxed
bulk disarrangement density%
\begin{equation}
\int_{\Omega }\left\vert (\nabla g-G)^{T}a\right\vert d\mathcal{L}^{N}
\label{final bulk energy (g,G) from u dot a}
\end{equation}%
for structured deformations $(g,G)$ of that region. The integral in \eqref{final bulk energy (g,G) from u dot a} represents the most economical way
of introducing jumps in the direction of $a$ while approaching in the limit
the given structured deformation $(g,G)$, including both jumps normal and
tangential to the discontinuity surfaces of approximating deformations $u$.
\ \ \ 

We note also the formula%
\begin{equation}
\max_{i=1,\ldots ,N}H(A,B,\delta _{i})=\left\Vert B-A\right\Vert _{\mathrm{%
row}\max }  \label{rowmax norm}
\end{equation}%
where on the left $\delta _{1},\ldots ,\delta _{N}$ denotes the standard
basis of $\mathbb{R}^{N}$ and on the right $\left\Vert B-A\right\Vert _{%
\mathrm{row}\max }$ denotes the maximum of the Euclidean norms of the rows
of $B-A$. The mapping $\left\Vert \cdot \right\Vert _{\mathrm{row}\max }:%
\mathbb{R}^{N\times N}\longrightarrow \mathbb{R}$ turns out to be a norm on $%
\mathbb{R}^{N\times N}$, and our interpretation of the integral in \eqref{final bulk energy (g,G) from u dot a} leads us to interpret the integral%
\begin{equation*}
\int_{\Omega }\left\Vert (\nabla g-G)(x)\right\Vert _{\mathrm{row}\max }d%
\mathcal{L}^{N}(x)
\end{equation*}%
as a bulk measure of disarrangements that takes into account at each $x\in
\Omega $ the direction $\delta _{i(x)}$ that maximizes the relaxed bulk
energy densities $H(\nabla g(x),G(x),\delta _{i})$ for $i=1,\ldots ,N$.
The bulk disarrangement density $\max_{i=1,\ldots ,N}H(A,B,\delta
_{i})=\left\Vert B-A\right\Vert _{\mathrm{row}\max }$ satisfies%
\begin{equation*}
\begin{split}
\max_{i=1,\ldots ,N}H(A,B,\delta _{i}) \leq  \inf \bigg\{ & \max_{i=1,\ldots ,N}\int_{J(u)}\left\vert [u](x)\cdot
\delta _{i}\right\vert d\mathcal{H}^{N-1}(x)\, : \\
&\; u\in SBV(Q;\mathbb{R}^{N}), u|_{\partial Q}=0\text{, }\nabla u=B-A\text{ }a.e.\bigg\},
\end{split}
\end{equation*}%
and need not be the relaxed bulk energy density corresponding to the initial
interfacial energy $ \displaystyle\max_{i=1,\ldots ,N}\int_{J(u)}\left\vert [u](x)\cdot
\delta _{i}\right\vert d\mathcal{H}^{N-1}(x)\,$. \

\bigskip
\textbf{Acknowledgments.} 
The authors warmly thank the CNA (NSF Grants No.\@ DMS-0405343 and DMS-0635983) at Carnegie Mellon University, Pittsburgh, USA and CAMGSD (FCT grant UID/MAT/04459/2013) at Instituto Superior Técnico, Lisbon, Portugal where this research was carried out.
The research of A.C.B, J.M., and M.M.\@ was partially supported by the Funda\c{c}\~{a}o para a Ci\^{e}ncia e a Tecnologia (Portuguese Foundation for Science and Technology) through the CMU-Portugal Program under grant FCT-UTA\_CMU/MAT/0005/2009 ``Thin Structures, Homogenization, and Multiphase Problems''.
The research of A.C.B.\@ was partially supported by the Funda\c{c}\~{a}o para a Ci\^{e}ncia e a Tecnologia through grant PEst\_OE/MAT/UI0209/2013.
The research of M.M. was partially supported by the European Research Council through the ERC Advanced Grant ``QuaDynEvoPro'', grant agreement no.\@ 290888. M.M. is a member of the Progetto di Ricerca GNAMPA-INdAM 2015 ``Fenomeni critici nella meccanica dei materiali: un approccio variazionale'' (INdAM-GNAMPA Project 2015 ``Critical phenomena in the mechanics of materials: a variational approach'').

\end{document}